# Reducing Congestion-Induced Renewable Curtailment with Corrective Network Reconfiguration in Day-Ahead Scheduling

Arun Venkatesh Ramesh
*Student Member, IEEE*
Department of Electrical and Computer Engineering
University of Houston
Houston, TX, USA
aramesh4@uh.edu

Xingpeng Li
*Member, IEEE*
Department of Electrical and Computer Engineering
University of Houston
Houston, TX, USA
xli82@uh.edu

*Abstract*— Renewable energy sources (RES) has gained a lot of interest recently. The limited transmission capacity serving RES often leads to network congestion since they are located in remote favorable locations. As a result, when poorly scheduled, the intermittent nature of RES may result in high curtailments of the free resource. Currently, grid operators utilize a static network when performing day-ahead scheduling and ignore transmission flexibility. This paper explores the possibility of utilizing network reconfiguration as a corrective action to reduce the transmission congestion and thereby the reduction of RES curtailments in day-ahead scheduling. To facilitate the RES integration in the grid, a stochastic *N-1* security-constrained unit-commitment with corrective network reconfiguration (SSCUC-CNR) is modelled. SSCUC-CNR model is studied on a modified IEEE 24-bus system with RES. The simulation results demonstrates that CNR not only leads to a lower cost solution by reducing network congestion but also facilitates RES integration by reducing congestion-induced curtailments in high penetration cases. Emission studies demonstrate that more green generators are committed resulting in reduced carbon emissions when CNR is implemented.

*Index Terms*—Corrective network reconfiguration, Flexible transmission, Renewable energy sources, Renewable curtailment, Stochastic programming, Post-contingency congestion relief.

## Nomenclature

| | |
|---|---|
| $g(n)$ | Set of generators connected to bus $n$. |
| $w(n)$ | Set of RES units connected to bus $n$. |
| $\delta^+(n)$ | Set of lines with bus $n$ as receiving bus. |
| $\delta^-(n)$ | Set of lines with bus $n$ as sending bus. |
| $UT_g$ | Minimum up time for generator $g$. |
| $DT_g$ | Minimum down time for generator $g$. |
| $c_g$ | Linear cost for generator $g$. |
| $c_g^{NL}$ | No-load cost for generator $g$. |
| $c_g^{SU}$ | Start-up cost for generator $g$. |
| $c_w^{pen}$ | Penalty for energy curtailed for RES $w$. |
| $P_g^{min}$ | Minimum output limit of generator $g$. |
| $P_g^{max}$ | Maximum output limit of generator $g$. |
| $P_{w,s}^{max}$ | Maximum capacity of RES $w$ in scenario $s$. |
| $R_g^{hr}$ | Regular hourly ramping limit of generator $g$. |
| $R_g^{SU}$ | Start-up ramping limit of generator $g$. |
| $R_g^{SD}$ | Shut-down ramping limit of generator $g$. |
| $R_g^{10}$ | 10-minute outage ramping limit of generator $g$. |
| $P_k^{max}$ | Long-term thermal line limit for line $k$. |
| $b_k$ | Susceptance of line $k$. |
| $P_k^{emax}$ | Emergency thermal line limit for line $k$. |
| $M$ | Real number with huge value |
| $\pi_s$ | Known probability of RES scenario $s$. |
| $P_{g,t,s}$ | Output of generator $g$ in time period $t$ and scenario $s$. |
| $P_{w,t,s}$ | RES $w$ output in time period $t$ and scenario $s$. |
| $u_{g,t}$ | Commitment status of generator $g$ in time period $t$. |
| $v_{g,t}$ | Start-up variable of generator $g$ in time period $t$. |
| $r_{g,t,s}$ | Reserve from generator $g$ in time period $t$. |
| $P_{k,t,s}$ | Line flow of line $k$ in time period $t$ and scenario $s$. |
| $\theta_{n,t,s}$ | Phase angle of bus $n$ in time period $t$ and scenario $s$. |
| $\theta_{m,t,s}$ | Phase angle of bus $m$ in time period $t$ and scenario $s$. |
| $d_{n,t}$ | Predicted demand of bus $n$ in time period $t$. |
| $P_{g,c,t,s}$ | Output of generator $g$ in time period $t$ and scenario $s$ after outage of line $c$ |
| $P_{w,c,t,s}$ | Output of RES $w$ in time period $t$ and scenario $s$ after outage of line $c$ |
| $P_{k,c,t,s}$ | Line flow of line $k$ in time period $t$ and scenario $s$ after outage of line $c$ |
| $\theta_{n,c,t,s}$ | Phase angle of bus $n$ in time period $t$ and scenario $s$ after outage of line $c$. |
| $\theta_{n,c,t,s}$ | Phase angle of bus $n$ in time period $t$ and scenario $s$ after outage of line $c$. |
| $z_{c,t,s}^k$ | Line status variable of line $k$ after outage of line $c$ in time period $t$. |

## I. Introduction

The importance on climate change and global warming in recent years has increased the investments in renewable sources of energy. The Paris climate deal set ambitious goals to reduce the carbon emissions by 2030 to limit the rise in global temperature [1]. Such directives place an emphasis on renewable energy sources (RES) as opposed to conventional fossil fuel plants. Typically, an increase in wind and solar generation is seen as favorable. However, the intermittent nature of RES due to weather brings challenges to the efficient and reliable grid operations [2].

During high penetration of RES, a flexible power system facilitates the integration of intermittent RES. This entails the usage of storage devices and flexible demand. Moreover, the requirement of favorable location and land implies RES are placed in remote locations. Therefore, even with the introduction of large-scale storage devices, the high penetration of RES results in curtailment due to network congestion. As a result, local generating sources utilizing fossil fuels are more utilized at the cost of RES curtailment. An effective smart grid and new technologies such as energy storage or flexible AC transmission System (FACTS) are required to utilize RES concurrently without spilling free RES.

To relieve network congestion and reduce RES curtailment, it requires transmission expansion planning to increase the transfer capability [3]. Another option is to redirect the power flow on the lines. This can be implemented through modifying the line parameters using FACTS devices [4] or network reconfiguration (NR) [5]-[7]. However, flexibility through expansion planning, energy sources and FACTS devices require expensive investment and maintenance. Therefore, the usage of NR is attractive to utilize the power produced by RES to meet the demand concurrently as it does not require any investment.

The grid network is built with redundancy to handle increasing future demand and maintain system reliability. This adds flexibility in the transmission network that is not fully considered. Currently, ISOs do not implement a dynamic network in day-ahead or real-time operations. This implies that the flexibility in power systems is provided by committing extra generators. In day-ahead operations, the ISOs utilize security-constrained unit commitment (SCUC) to commit generators with a goal to minimize operational costs while respecting physical and reliability constraints. Thus, facilitating network reconfiguration (NR) can increase the system flexibility while reducing overall cost.

NR can be of both preventive action and corrective action. However, concerns that NR causes a big network disturbance, stability issues and circuit breaker degradation makes corrective network reconfiguration (CNR) more attractive as it is only implemented after a contingency has occurred. Prior research shows the cost-saving benefits of NR [8] and CNR [9] due to the increased feasible set of solutions for the SCUC problem.

CNR first introduced in [10] is attractive in reducing line overloading and relieving congestion [9], [11]. CNR implementation is also scalable to large-scale networks in real-time operations by using practical and innovative heuristic methods for post-contingency violation reduction [12] and post-contingency network congestion management [13]-[14]. In addition, CNR offers increased network flexibility as shown in [15] where it was implemented on an industry case using an in-house industry software. The impact of NR on high penetrative wind models were studied in [16]-[19]. [20] provides a real-time implementation of enhancing optimal power flow by incorporating CNR in economic dispatch to facilitate integration of RES in the grid. However, the effect of SCUC with CNR on high penetrative RES network and RES curtailment studies has not been performed. Due to the high variability of RES, it requires solution which is satisfied in multiple scenarios. Therefore, a stochastic implementation through a known probability distribution of multiple scenarios is considered for a feasible solution as seen in [21]-[24].

In [25], optimal NR is implemented through a bi-level stochastic implementation to solve large scale networks. However, this paper does not consider the use of reconfiguration as a corrective action and post-contingency constraints were not modelled. Other viable technologies for reducing RES curtailments is through FACTS to reduce network congestion [26] and the use of energy storage [27].

Hence, in this paper, we study the benefits of CNR in SCUC model on RES curtailment by incorporating a stochastic model with multiple scenarios with a known probability distribution. Though the test case considered in this paper is mainly related wind energy, this work can also be implemented for other variable RES such as solar energy. The rest of this paper is organized as follows. Section II depicts the stochastic model of day-ahead *N*-1 SCUC with and without CNR for a system with renewable generation. Section III details the RES scenarios and the data used. A discussion of results is shown in Section IV and conclusions drawn from the results are summarized in Section V.

## II. MATHEMATICAL MODEL

This paper proposes two models to determine the least-cost generator commitment and dispatch solutions in day-ahead scheduling considering multiple system scenarios: a stochastic-SCUC (SSCUC) model and a stochastic-SCUC with CNR (SSCUC-CNR) model. SSCUC and SSCUC-CNR are based on the simplified DC power flow model and they are subject to base-case and post-contingency physical requirements of the traditional generators, renewable generation limits and transmission constraints while meeting the demand. The utilization of free RES output is directly related to reducing the cost and hence the curtailments are lower in base-case. But, both SSCUC and SSCUC-CNR solution leads to high post-contingency RES curtailment as it is not considered in the objective. Since the study is focused on reducing or eliminating RES curtailments, a penalty cost, $c_w^{pen}$, was added for post-contingency curtailment as shown in (1).

$$Min: \sum_{g,t}(c_g^{NL} u_{g,t} + c_g^{SU} v_{g,t} + \sum_s(\pi_s c_g P_{g,t,s})) + \sum_{w,c,t,s}(\pi_s c_w^{pen}(P_w^{max} - P_{w,c,t,s})) \quad (1)$$

The base-case generation constraints, (2)-(13), consist of the min-max limits of generator output, reserve limits, generator ramping requirements, minimum up-down time, generator start-up and commitment constraints bounded by integrality and finally the maximum RES generation constraints.

$$P_g^{min} u_{g,t} \leq P_{g,t,s}, \forall g,t,s \quad (2)$$
$$P_{g,t,s} + r_{g,t,s} \leq P_g^{max} u_{g,t}, \forall g,t,s \quad (3)$$
$$0 \leq r_{g,t,s} \leq R_g^{10} u_{g,t}, \forall g,t,s \quad (4)$$
$$\sum_{q \in G} r_{q,t,s} \geq P_{g,t,s} + r_{g,t,s}, \forall g,t,s \quad (5)$$
$$P_{g,t,s} - P_{g,t-1,s} \leq R_g^{hr} u_{g,t-1} + R_g^{SU} v_{g,t}, \forall g,t,s \quad (6)$$
$$P_{g,t-1,s} - P_{g,t,s} \leq R_g^{hr} u_{g,t} + R_g^{SD}(v_{g,t} - u_{g,t} + u_{g,t-1}), \forall g,t,s \quad (7)$$
$$\sum_{q=t-UT_g+1}^{t} v_{g,q} \leq u_{g,t}, \forall g, t \geq UT_g \quad (8)$$
$$\sum_{q=t+1}^{t+DT_g} v_{g,q} \leq 1 - u_{g,t}, \forall g, t \leq T - DT_g \quad (9)$$
$$v_{g,t} \geq u_{g,t} - u_{g,t-1}, \forall g,t \quad (10)$$
$$v_{g,t} \in \{0,1\}, \forall g,t \quad (11)$$
$$u_{g,t} \in \{0,1\}, \forall g,t \quad (12)$$
$$0 \leq P_{w,t,s} \leq P_{w,s}^{max}, \forall w,t,s \quad (13)$$

The base-case transmission constraints, (14)-(16), consist of DC power flow equation, the min-max line long-term thermal limits, and the nodal power balance equations with renewable generation injection.

$$P_{k,t,s} - b_k(\theta_{n,t,s} - \theta_{m,t,s}) = 0, \forall k,t,s \quad (14)$$
$$-P_k^{max} \leq P_{k,t,s} \leq P_k^{max}, \forall k,t,s \quad (15)$$
$$\sum_{g \in g(n)} P_{g,t,s} + \sum_{k \in \delta^+(n)} P_{k,t,s} - \sum_{k \in \delta^-(n)} P_{k,t,s} = d_{n,t} - \sum_{w \in w(n)} P_{w,t,s}, \forall n,t,s \quad (16)$$

The post-contingency case generator constraints, (17)-(21), models the generator 10-minute ramp up-down and min-max limits and renewable generation limit after the outage of line $c$.

$$P_{g,t,s} - P_{g,c,t,s} \leq R_g^{10} u_{g,t}, \forall g,c,t,s \quad (17)$$
$$P_{g,c,t,s} - P_{g,t,s} \leq R_g^{10} u_{g,t}, \forall g,c,t,s \quad (18)$$
$$P_g^{min} u_{g,t} \leq P_{g,c,t,s}, \forall g,c,t,s \quad (19)$$
$$P_{g,c,t,s} \leq P_g^{max} u_{g,t}, \forall g,c,t,s \quad (20)$$
$$0 \leq P_{w,c,t,s} \leq P_{w,s}^{max}, \forall w,t,s \quad (21)$$

The post-contingency case transmission constraints, (22)-(27), where the nodal balance is maintained under a line outage through (22). Branch power flow equations and limits without CNR is modelled in (23)-(24) whereas the branch power flow equations and limits with CNR is modelled in (25)-(27). The linearity of post-contingency power flow equations, (24)-(25), are maintained by the 'big-M' method. The binary decision variable, $z_{c,t}^k$, represents the CNR action where the value 0 represents line is disconnected from the system and the value of 1 indicates line is available. These contingencies are modelled for all non-radial lines. A restriction on the number of CNR actions in each post-contingency case is introduced through (28) to reduce system disturbance.

$$\sum_{g \in g(n)} P_{g,c,t,s} + \sum_{k \in \delta^+(n)} P_{k,c,t,s} - \sum_{k \in \delta^-(n)} P_{k,c,t,s} = d_{n,t} - \sum_{w \in w(n)} P_{w,c,t,s}, \forall n,c,t,s \quad (22)$$
$$P_{k,c,t,s} - b_k(\theta_{n,c,t,s} - \theta_{m,c,t,s}) = 0, \forall k,c,t,s \quad (23)$$
$$-P_k^{emax} \leq P_{k,c,t,s} \leq P_k^{emax}, \forall k,c,t,s \quad (24)$$
$$P_{k,c,t,s} - b_k(\theta_{n,c,t,s} - \theta_{m,c,t,s}) + (1-z_{c,t,s}^k)M \geq 0, \forall k,c,t,s \quad (25)$$
$$P_{k,c,t,s} - b_k(\theta_{n,c,t,s} - \theta_{m,c,t,s}) - (1-z_{c,t,s}^k)M \leq 0, \forall k,c,t,s \quad (26)$$
$$-P_k^{emax} z_{c,t,s}^k \leq P_{k,c,t,s} \leq z_{c,t,s}^k P_k^{emax}, \forall k,c,t,s \quad (27)$$
$$\sum_k (1-z_{c,t,s}^k) \leq 1, \forall k,c,t,s \quad (28)$$

SSCUC is modelled by (1)-(24) whereas the SSCUC-CNR is represented by (1)-(22) and (25)-(28). The effectiveness of CNR is demonstrated considering only line outages. It can be noted that the commitment schedule and the start-up variable of the generators are the same across all scenarios, $\forall s$, therefore the constraints (8)-(12) are not scenario based.

### III. TEST CASE: IEEE 24-BUS SYSTEM WITH RES

To study the effect of CNR, the IEEE 24-bus network is utilized for testing [28]. The base system contains 24 buses, 33 generators and 38 branches. However, the system was modified to include three wind farms located at bus 12, bus 16 and bus 22 to study the effect of network constraints on RES curtailment. Five different scenarios are considered for wind generation; and the base total system renewable generation over 24 hours for each scenario are represented in Fig. 1.

The system-wide RES output for various penetration level is represented in Fig. 2. The base total RES output was modified to obtain five cases considered for the study and can be classified using the peak load period penetration as ~15%, ~30%, ~50%, ~60 and ~80%. Apart from the wind generation, the total generation capacity from traditional units is 3,393 MW and the system peak load is 2,270 MW. The wind output was assumed to be constant for each three-hour-period due to the computational complexity of CNR for this study.

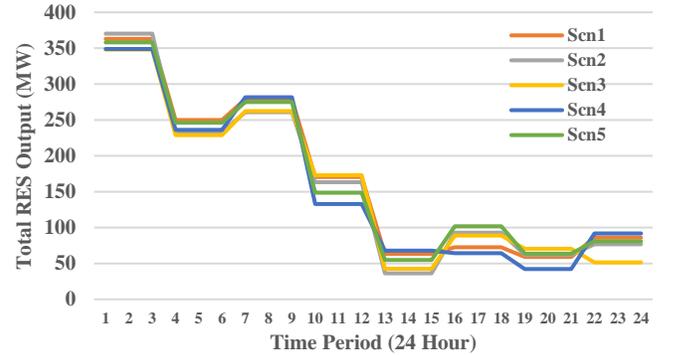
Fig. 1. The base total RES output for each scenario.

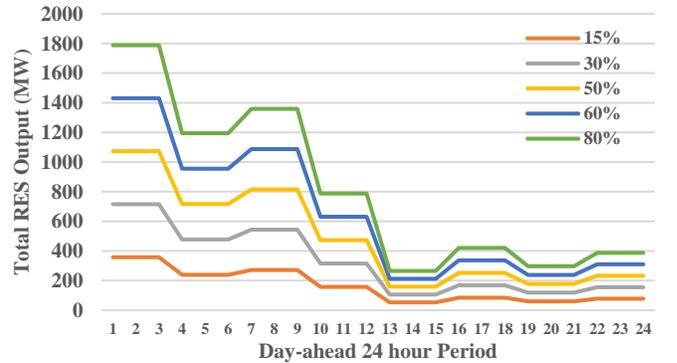
Fig. 2. System-wide RES generation for various penetration levels.

### IV. RESULTS AND ANALYSIS

The mathematical model is implemented using AMPL and solved using Gurobi with a MIPGAP of 0.01 for a 24-hour (day-ahead) load period. The computer with Intel® Xeon(R) W-2195 CPU @ 2.30GHz, and 128 GB of RAM was utilized.

#### A. Rationale for Penalty Cost:

Initially, the system with peak penetration of 30% was studied under two cases, complete wind output usage (CWOU) that uses up all available wind power and variable wind output usage (VWOU) that allows the system to curtail some wind power for both models with and without CNR. Table I and Table II, shows the total costs, the base-case curtailments (BCC) and the expected post-contingency curtailments (PCC) for day-ahead SSCUC and SSCUC-CNR respectively. The BCC is aggregated over all periods, $\forall t$ and RES units, $\forall w$. Similarly, the PCC is aggregated over all periods, $\forall t$, and RES units, $\forall w$, and then averaged over all contingencies, $\forall c$. Since, this is a multi-scenario stochastic implementation, the probability of scenarios is utilized to obtain BCC and PCC as shown in (29) and (30), respectively.

$$BCC = \left(\sum_{w,t,s}(\pi_s(P_w^{max} - P_{w,t,s}))\right) \quad (29)$$
$$PCC = \left(\sum_{w,c,t,s}(\pi_s(P_w^{max} - P_{w,c,t,s}))\right)/n_c \quad (30)$$

From the initial assessment, SSCUC-CNR offers lower total costs. But from VWOU, both SSCUC and SSCUC-CNR are susceptible to heavy renewable curtailment since there is no cost associated with PCC in the objective when minimizing operational costs. It is also seen that the total cost for CWOU and VWOU is the same for the respective implementations. Therefore, introducing a penalty cost to limit PCC eliminates the curtailment in post-contingency cases without increasing total cost. Further studies in this section leverage the penalty for PCC.

TABLE I. PENALTY COST STUDIES FOR SSCUC

|  | No PCC penalty | | With PCC Penalty |
| --- | --- | --- | --- |
|  | CWOU | VWOU |  |
| Total cost ($) | 201,769 | 201,769 | 201,769 |
| BCC (MW) | NA | 0 | 0 |
| PCC (MW) | NA | 9.14 | 0 |

TABLE II. PENALTY COST STUDIES FOR SSCUC-CNR

|  | No PCC penalty | | With PCC Penalty |
| --- | --- | --- | --- |
|  | CWOU | VWOU |  |
| Total cost ($) | 177,170 | 177,170 | 177,170 |
| BCC (MW) | NA | 0 | 0 |
| PCC (MW) | NA | 80.90 | 0 |

### B. RES Penetration Sensitivity Studies:

The renewable energy penetration based sensitivity studies are performed, and the associated results are presented in Fig. 3 and Fig. 4. As shown in Fig. 3, the general trend observed is that (i) the total cost reduces as more free renewable power was utilized and (ii) SSCUC-CNR always offered lower cost solutions compared to SSCUC. This demonstrates that CNR alleviates network congestion and reduce congestion-induced cost by increasing the transmission flexibility.

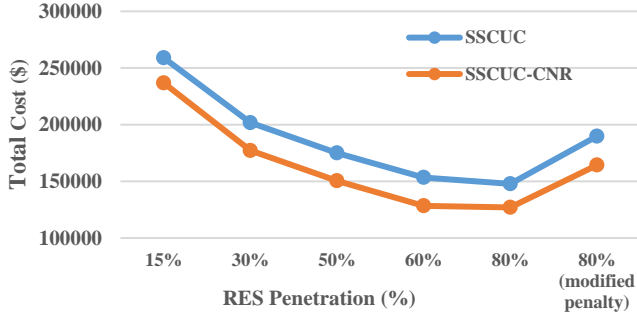

Fig. 3. Total Cost in $ under various penetration levels.

As shown in Fig. 4, the network utilizes all available renewable power in low penetration levels. However, RES curtailments are observed for both base case and contingency cases when RES penetration level is above 30%; it is also observed that CNR actions alleviate the post-contingency congestions for high penetration levels, 50%-80%. However, under very-high penetration of 80%, CNR alone is not beneficial as both PCC and BCC are higher with SSCUC-CNR against traditional SSCUC. This can be characterized by the cost saving offered through congestion alleviation that provides a much lower cost (even it includes the penalty for curtailment). Therefore, increasing the penalty costs for 80% RES penetration resulted in a reduction of PCC from 323 MW to 153 MW and a reduction of BCC from 1047 MW to 655 MW for SSCUC-CNR. However, there are no changes for SSCUC; instead, the total cost increases marginally as a result of higher penalty. Hence, a combination of dynamic penalty factor along with CNR may be more beneficial.

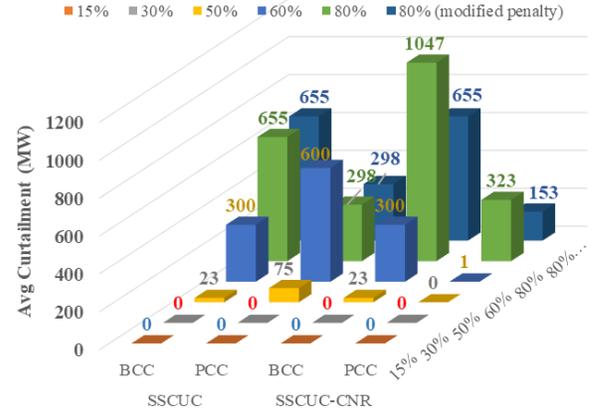

Fig. 4. RES curtailment under various penetration levels.

### C. Carbon Emission Reduction:

One of the key aspects of integrating renewables in the system is the reduction of emissions. The emission data for the generators are used to highlight the reduced emissions. The base-case generation outputs were used to calculate the total net heat and emission of each generator for the test system [28]. When averaged over multiple scenarios, it is seen from Fig. 5 that SSCUC-CNR leads to significantly lower carbon emissions at high penetration, 60%-80%, of RES. In comparison, SSCUC shows an increase in emissions at high penetration of RES due to higher curtailments. This implies more traditional generators are used to meet the demand thereby increasing carbon emissions.

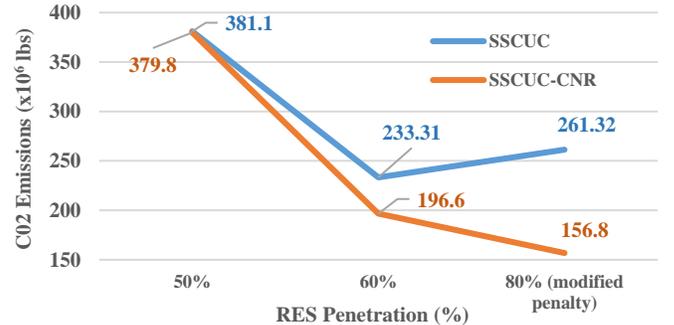

Fig. 5. System carbon emission under various penetration levels.

At 50% peak penetration, both SSCUC and SSCUC-CNR leads to similar BCC; however, the scheduling solutions with SSCUC-CNR resulted in $1.3 \times 10^6$ lbs reduction of carbon emissions as compared to SSCUC. This implies that CNR is also beneficial at lower penetration level of RES as it leads to the utilization of efficient low-cost "green" generators through the alleviation of network congestion.

### D. Corrective Network Reconfiguration Strategy:

The peak RES penetration of 80% shows high PCC in periods 1-3 and 7-9 due to the load profile, intermittent nature of RES and network congestion. However, periods 1-3 shows the most reconfiguration action taken. In total there were 64 line outage cases across 5 scenarios, 38 lines in periods 1-3 that required

CNR action. Since only one line is removed as an action, a pattern to note here is that line 5 [bus 2 - bus 6], line 16 [bus 10 - bus 11], and line 30 [bus 17 - bus 18] were common choices to remove from the network.

On closer observation, firstly, these lines are closer to where RES are located in the network. Secondly, the bottleneck lines are typically line 10 [bus 6 - bus 10] and line 23 [bus 14 - bus 16]. The above CNR actions help relieve the congested lines which in turn reduces RES curtailments.

## V. Conclusions

The increase in RES in the network is key to addressing climate change issues. Due to the intermittent nature of RES, smart grids are required to facilitate the integration of RES. To avoid undesired congestion-induced curtailment of free energy, a flexible network is required. The day-ahead operational procedure still uses a static network which impedes further deployment of renewables in the grid. Hence, to reduce the RES curtailments of available renewable generation while considering the intermittent nature of RES, we proposed a corrective dynamic network for contingencies implemented in SSCUC-CNR in this paper.

It was observed that SSCUC-CNR provides more transfer capability of the network thereby avoiding congestion-induced and contingency-induced RES curtailment in high-penetration of RES. Along with reduction of curtailment, SSCUC-CNR also lowers the cost of operation, and reduces green-house gas emissions. Numerical simulations also showed that SSCUC-CNR is also beneficial in moderate-penetration of RES by committing efficient green (less emission) generators which reduces the overall carbon emissions in a day-ahead schedule.

The future work to be considered includes the scalability of this model for large power system networks. It is known that the addition of CNR leads to increased solution complexity which can be addressed by decomposing the problem as a master-slave problem to reduce the computational burden. Another aspect being considered is modelling energy storage to reduce RES curtailments.

## VI. References

[1] Paris Agreement, 2030 Climate and Energy Framework, [Online]. Available: https://ec.europa.eu/clima/policies/strategies/2030
[2] M. Milligan et al., "Wind power myths debunked," IEEE Power Energy Mag., vol. 7, no. 6, pp. 89–99, Nov./Dec. 2009.
[3] J. C. Villumsen, G. Brønmo and A. B. Philpott, "Line capacity expansion and transmission switching in power systems with large-scale wind power," *IEEE Transactions on Power Systems*, vol. 28, no. 2, pp. 731-739, May 2013.
[4] A. Nasri, A. J. Conejo, S. J. Kazempour, and M. Ghandhari, "Minimzing wind power splillage using an OPF with FACTS devices," IEEE Trans. Power Syst., vol. 29, no. 5, pp. 2150–2159, Sep. 2014.
[5] E. B. Fisher, R. P. O'Neill and M. C. Ferris, "Optimal Transmission Switching," *IEEE Transactions on Power Systems*, vol. 23, no. 3, pp. 1346-1355, Aug. 2008.
[6] K. W. Hedman, S. S. Oren and R. P. O'Neill, "A review of transmission switching and network topology optimization," *2011 IEEE Power and Energy Society General Meeting*, San Diego, CA, 2011, pp. 1-7.
[7] H. Glavitsch, "State of the art review: switching as means of control in the power system", *INTL. JNL. Elect. Power Energy Syst.*, vol. 7, no. 2, pp. 92-100, Apr 1985.
[8] K. W. Hedman, R. P. O'Neill, E. B. Fisher and S. S. Oren, "Optimal Transmission Switching with Contingency Analysis," *IEEE Transactions on Power Systems*, vol. 24, no. 3, pp. 1577-1586, Aug. 2009.
[9] Arun Venkatesh Ramesh, Xingpeng Li, "Security constrained unit commitment with corrective transmission switching", *North American Power Symposium*, Wichita, KS, USA, October 2019.
[10] Koglin, and Müller. "Corrective Switching: A New Dimension in Optimal Load Flow." *International Journal of Electrical Power and Energy Systems*, 4.2 (1982): 142-49.
[11] Wei Shao and V. Vittal, "Corrective switching algorithm for relieving overloads and voltage violations," *IEEE Transactions on Power Systems*, vol. 20, no. 4, pp. 1877-1885, Nov. 2005.
[12] Xingpeng Li, P. Balasubramanian, M. Sahraei-Ardakani, M. Abdi-Khorsand, K. W. Hedman and R. Podmore, "Real-Time Contingency Analysis with Corrective Transmission Switching," *IEEE Transactions on Power Systems*, vol. 32, no. 4, pp. 2604-2617, July 2017.
[13] Xingpeng Li and Kory W. Hedman, "Enhanced Energy Management System with Corrective Transmission Switching Strategy— Part I: Methodology," *IEEE Transactions on Power Systems*, vol. 34, no. 6, pp. 4490-4502, Nov. 2019.
[14] Xingpeng Li and Kory W. Hedman, "Enhanced Energy Management System with Corrective Transmission Switching Strategy— Part II: Results and Discussion," *IEEE Transactions on Power Systems*, vol. 34, no. 6, pp. 4503-4513, Nov. 2019.
[15] J. D. Lyon, S. Maslennikov, M. Sahraei-Ardakani, Tongxin Zheng, E. Litvinov, Xingpeng Li, P. Balasubramanian, And K. W. Hedman, "Harnessing Flexible Transmission: Corrective Transmission Switching for ISO-NE," *IEEE Power and Energy Technology Systems Journal*, vol. 3, no. 3, pp. 109-118, Sept. 2016.
[16] J. C. Villumsen, G. Brønmo and A. B. Philpott, "Line capacity expansion and transmission switching in power systems with large-scale wind power," *IEEE Transactions on Power Systems*, vol. 28, no. 2, pp. 731-739, May 2013.
[17] T. Lan and G. M. Huang, "Transmission line switching in power system planning with large scale renewable energy," *2015 First Workshop on Smart Grid and Renewable Energy (SGRE)*, Doha, 2015, pp. 1-6.
[18] Nikoobakht, A., Mardaneh, M., Aghaei, J., Guerrero-Mestre, V., Contreras, J., & Nikoobakht, A. (2017). Flexible power system operation accommodating uncertain wind power generation using transmission topology control: an improved linearised AC SCUC model. *IET Generation, Transmission and Distribution*, *11*(1), 142–153.
[19] F. Qiu and J. Wang, "Chance-Constrained Transmission Switching with Guaranteed Wind Power Utilization," in *IEEE Transactions on Power Systems*, vol. 30, no. 3, pp. 1270-1278, May 2015.
[20] Xingpeng Li and Qianxue Xia, "Stochastic Optimal Power Flow with Network Reconfiguration: Congestion Management and Facilitating Grid Integration of Renewables", *IEEE PES T&D Conference & Exposition*, Chicago, IL, USA, Apr. 2020, accepted for publication.
[21] J. Hetzer, D. C. Yu and K. Bhattarai, "An Economic Dispatch Model Incorporating Wind Power," *IEEE Transactions on Energy Conversion*, vol. 23, no. 2, pp. 603-611, June 2008.
[22] M. Vrakopoulou, K. Margellos, J. Lygeros and G. Andersson, "A Probabilistic Framework for Reserve Scheduling N-1 Security Assessment of Systems With High Wind Power Penetration," *IEEE Transactions on Power Systems*, vol. 28, no. 4, pp. 3885-3896, Nov. 2013.
[23] A. Nikoobakht, J. Aghaei and M. Mardaneh, "Optimal transmission switching in the stochastic linearised SCUC for uncertainty management of the wind power generation and equipment failures," in *IET Generation, Transmission & Distribution*, vol. 12, no. 17, pp. 4060-4060, 30 9 2018.
[24] O. Ziaee and F. Choobineh, "Optimal Location-Allocation of TCSCs and Transmission Switch Placement Under High Penetration of Wind Power," *IEEE Transactions on Power Systems*, vol. 32, no. 4, pp. 3006-3014, July 2017.
[25] J. Shi and S. S. Oren, "Stochastic Unit Commitment With Topology Control Recourse for Power Systems With Large-Scale Renewable Integration," *IEEE Transactions on Power Systems*, vol. 33, no. 3, pp. 3315-3324, May 2018.
[26] Sang and M. Sahraei-Ardakani, "Stochastic Transmission Impedance Control for Enhanced Wind Energy Integration," *IEEE Transactions on Sustainable Energy*, vol. 9, no. 3, pp. 1108-1117, Jul. 2018.
[27] X. Dui, G. Zhu and L. Yao, "Two-Stage Optimization of Battery Energy Storage Capacity to Decrease Wind Power Curtailment in Grid-Connected Wind Farms," *IEEE Transactions on Power Systems*, vol. 33, no. 3, pp. 3296-3305, May 2018.
[28] C. Grigg *et al.*, "The IEEE Reliability Test System-1996. A report prepared by the Reliability Test System Task Force of the Application of Probability Methods Subcommittee," *IEEE Transactions on Power Systems*, vol. 14, no. 3, pp. 1010-1020, Aug. 1999.